\documentclass[11pt]{article}
\usepackage{amssymb}
\usepackage{amsfonts}
\usepackage{amsmath}
\usepackage{amsthm}
\usepackage{enumerate}

\usepackage{srcltx}
\newtheorem{theorem}{Theorem}
\newtheorem{lemma}{Lemma}
\newtheorem{corollary}{Corollary}
\newtheorem{example}{Example}
\newtheorem{remark}{Remark}
\newcommand\cA{{\cal A}}
\newcommand\cC{{\cal C}}

\newcommand\cL{{\cal L}}
\newcommand\cB{{\cal B}}
\newcommand\cP{{\cal P}}
\newcommand\cQ{{\cal Q}}

\newcommand\cN{{\cal N}}
\newcommand\cM{{\cal M}}
\newcommand\cX{{\cal X}}
\newcommand\cY{{\cal Y}}
\newcommand\cD{{\cal D}}
\newcommand\cS{{\cal S}}
\newcommand\cR{{\cal R}}

\newcommand{\wh}{\widehat}
\newcommand{\wt}{\widetilde}

\def\bbr{{\mathbb R}}
\def\bbn{{\mathbb N}}

\def\bea{\begin{eqnarray}}
\def\eea{\end{eqnarray}}

\def\text#1{\hbox{#1}}
\def\proof{{\noindent \bf Proof. }}
\def\endproof{\mbox{\ $\Box$}}

\def\C{{\bf C}}

\def\E{{\bf E}}
\def\P{{\bf P}}

\def\R{{\bf R}}
\def\L{{\bf L}}
\def\Chi{{\bf 1}}
\def\d{\mathrm{d}}
\def\build #1_#2{\mathrel{\mathop{\kern 0pt #1}\limits_{#2}}}
\newcommand{\zs}[1]{{\mathchoice{#1}{#1}{\lower.25ex\hbox{$\scriptstyle#1$}}
{\lower0.25ex\hbox{$\scriptscriptstyle#1$}}}}

\def\proof{{\noindent \bf Proof. }}
\def\endproof{\mbox{\ $\qed$}}

\begin{document}
\title{General  model selection estimation of a periodic regression with a Gaussian noise}

\author{Victor Konev\thanks{
Department of Applied Mathematics 
and Cybernetics, Tomsk State University, Lenin street 36,
 634050 Tomsk, Russia. e-mail: vvkonev@mail.tsu.ru 
}
\and Serguei Pergamenchtchikov
\thanks{Laboratoire de Math\'ematiques Rapha\"el Salem,
UMR 6085 CNRS, Universit\'e de Rouen,
Avenue de l'Universit\'e, BP.12,
76801 Saint Etienne du Rouvray (France).\newline
email:  Serge.Pergamenchtchikov@univ-rouen.fr}
}
\date{}
\maketitle

\begin{abstract}
This paper considers the problem of estimating a periodic
function in a continuous time regression model with an additive
stationary gaussian noise having unknown correlation function.
 A general model selection procedure on
the basis of arbitrary projective estimates, which does not need
the knowledge of the noise correlation function,
 is proposed. A
non-asymptotic upper bound for $\cL_2$-risk (oracle inequality)
has been derived under mild conditions on the noise. For the 
Ornstein-Uhlenbeck noise the risk upper bound is shown to be uniform
in the nuisance parameter.
In the case of gaussian white noise the
constructed procedure has some advantages
 as compared with the
procedure based on the least squares estimates (LSE).
 The asymptotic minimaxity of the estimates
has been proved.
The proposed model selection scheme is extended also to
the estimation problem based on the discrete data applicably
to the situation when high frequency sampling can not be provided.
\end{abstract}



\bibliographystyle{plain}

\newpage

\section{Introduction}
Consider a regression model in continuous time
 \begin{equation}\label{1.1}
  \d y_t\,=\,S(t)\d t\,+\,\d \xi_t\,,
 \end{equation}
where $S(t)$ is an unknown $1$-periodic function in the space
$\cL_2[0,1]$, $(\xi_t)_\zs{t\ge 0}$ is a  continuous
gaussian process with zero mean and
such that for each $n\ge 1$
 the stochastic integral
$\int^n_0\,f(t)\d \xi_t$ is well-defined for any non-random
function $f$ from $\cL_2[0,n]$. The correlation function of noise
$\xi_t$ is unknown.

This process can be modeled in different ways.

\begin{example}\label{Ex.1.1}
 $\xi_\zs{t}$ is
a scalar non-explosive
Ornstein-Uhlenbeck process defined by the equation
\begin{equation}\label{1.2}
\d \xi_\zs{t}=\theta \xi_\zs{t}\d t+\d w_t\,,
\end{equation}
where $(w_\zs{t})_\zs{t\ge 0}$ is a standard brownian motion
and $\theta\le 0$ is unknown parameter;
 the initial value $\xi_\zs{0}\sim \cN(0,1/2|\theta|)$ if
$\theta< 0$ and $\xi_\zs{0}=0$ if $\theta=0$.
\end{example}

\begin{example}\label{Ex.1.2}
$\xi_\zs{t}$ is
a stationary
autoregressive
process of order $q\ge 2$
 satisfying the stochastic differential
equation
\begin{equation}\label{1.3}
\xi^{(q)}_\zs{t}=\theta_\zs{1} \xi^{(q-1)}_\zs{t}
+\ldots+\theta_\zs{q} \xi_\zs{t}+\dot{w}_t\,.
\end{equation}
Here $(\dot{w}_\zs{t})_\zs{t\ge 0}$ is  a white gaussian  noise
and
the unknown vector
$\theta=(\theta_\zs{1},\ldots,\theta_\zs{q})'$ belongs to stability
region of the process
\begin{equation}\label{1.4}
\cA = \{\theta\in\bbr^q\,:\, \max_\zs{1\le i\le q} \mbox{Re} \lambda_\zs{i}(\theta)
< 0 \}\,,
\end{equation}
where $(\lambda_\zs{i}(\theta))_\zs{1\le i \le q}$ are
eigenvalues of the matrix
\begin{equation}\label{1.5}
A=A(\theta) =
\left(
\begin{array}{ccc}
\theta_\zs{1}&\ldots&\theta_\zs{q}\\
&I_\zs{q-1} & 0
\end{array}
 \right)\,;
\end{equation}
$I_\zs{q}$ is the identity matrix of order $q$.
\end{example}

Models of type \eqref{1.1} and their discrete-time analogues have
been studied by a number of authors (see,  Efroimovich (1999),
Liptser and Shyraev (1974), Konev and Pergamenshchikov (2003),
Nemirovskii (2000) and references therein). The
estimation problem of periodic signal $S(t)$ in model
\eqref{1.1}--\eqref{1.2} has been thoroughly studied in the case,
when $(\xi_t)_\zs{t\ge 0}$ is a white gaussian noise (see, for
example, Ibragimov and Hasminskii (1981) for details and further references).

A discrete-time counterpart of model \eqref{1.1}--\eqref{1.2} was
applied in the econometrical problems for modeling the
consumption as a function of income Golfeld and Quandt (1972).

As is well known, the problem of
nonparametric estimation of $S(t)$ comprises the following three
statements: the function estimation
 at a fixed point $t_0$, estimation  in the
uniform metric and  in the integral metric. The first two problems
are usually solved by making use of the kernel and local
polynomial estimates. This paper focuses on the third setting with
the quadratic metric. The estimation in the integral metric is
based, as a rule, on the projective estimates which were first
proposed in Chenstov (1962) for estimating the distribution density in a
scheme of i.i.d. observations. The heart of this method is to
approximate the unknown function with a finite Fourier series.
Applying the projective estimates to the regression model \eqref{1.1}
with a white noise leads to the optimal convergence rate in
$\L_2(0,1)$ provided that the smoothness of $S$ is known (see for
example Ibragimov and Hasminskii (1981)). Another adaptive approach based on the model
selection method (see for example, Barron et al. (1999), Baraud (2000), Birg\'e and Massart (2001)
and  Fourdrinier and Pergamenshchikov (2007))
 enables one to study
this problem in the nonasymptotic setting
 when the smoothness of  function $S$ is
unknown. It should be noted that this method can be used also for
model \eqref{1.1} under the condition that the correlation
function $\E\xi_t\xi_s$ is exactly known and besides the unknown
function $S$ belongs to the subspace spanned by its eigenfunctions
(see, Theorem 1, p. 11 in Birg\'e and Massart (2001)). In our case, when the noise
correlation function is unknown, this method can not be applied.
 This paper develops a general model selection method
for the regression scheme \eqref{1.1} with unknown correlation properties.

Note that the usual nonasymptotic selection model procedure proposed in
 Barron et al. (1999), Baraud (2000), Birg\'e and Massart (2001)
is based on the least square estimators (LSE) which, as was shown
in  Goloubev (1982) and  Pinsker (1981), are not efficient in the problem of
nonparametric regression. Our approach is close to the general
model selection method proposed in Fourdrinier and Pergamenshchikov (2007) for discrete time
models with spherically symmetric errors which allows one to use any projective estimators in the model
selection procedure including the LSE. In Section 2 we propose a general
model selection procedure for a regression scheme in continuous time 
\eqref{1.1} with unknown correlation structure of the gaussian noise.
In Theorem~\ref{Th.2.1},
under some loose conditions on the noise,
 we
establish a nonasymptotic upper bound for the quadratic risk in
which the principal term is minimal over the set of all admissible
basic estimates. The inequalities of this type are usually called
{\sl oracle}.

In the case of the Ornstein-Uhlenbeck noise \eqref{1.2}, the risk upper bound
is shown to be uniform in the nuisance parameter (Corollary~\ref{Co.2.2}).

The rest of the paper is organized as follows. In Section 3 we
consider  case of white gaussian noise $\xi_t$ and show that the
possibility to choose different projective estimators in the
procedure may lead to a sharper upper bound for the mean square
estimation accuracy. In Section 4 the upper bound and the lower
bounds for the minimax quadratic risk are obtained under the
assumption that the smoothness of $S$ is unknown. In Section 6 we
consider the estimation problem for the regression model
\eqref{1.1} assuming that it is accessible for observations only at
discrete times $t_k=k/p$, $k=0,1,\ldots$. Such observation scheme
is more appropriate in a number of applications, where one can not
provide high frequency data sampling. 
Theorems~\ref{Th.6.1} establishs the nonasymptotic oracle inequalities in this
case. Appendix contains some technical
results.
\\[2mm]

\section{Nonasymptotic estimation}

\medskip

In this section we consider the estimation problem for the model
\eqref{1.1} in nonasymptotic setting, i.e. assuming that the
estimator of $S$ is based on the observations $(y_\zs{t})_\zs{0\le
t\le n}$ with a  fixed duration $n$. For this we apply the general
model selection approach proposed in Fourdrinier and Pergamenshchikov (2007) for the
discrete-time regression model.

 First we introduce some notations. Let
$\cX$ be the Hilbert space  of square integrable $1$-periodic functions on
$\bbr$ with the usual scalar product
$$
(x,y)\,=\,\int^1_0\,x(t)\,y(t)\,\d t
$$
and $(\phi_j)_\zs{j\ge 1}$ be a system of orthonormal functions in $\cX$,
i.e. $(\phi_i,\phi_j)=0$, if $i\ne j$ and $\|\phi_i\|^2=(\phi_i,\phi_i)=1$.

Then we impose the following additional conditions on the noise
$(\xi_t)_\zs{t\ge 0}$ in \eqref{1.1}. Assume that

\noindent $\C_1)$  {\sl For each $n\ge 1$ and $k\ge 1$ the vector
$\zeta(n)=(\zeta_\zs{1}(n),\ldots,\zeta_\zs{k}(n))^\prime$ with
components
 \begin{equation}\label{2.1}
  \zeta_\zs{j}(n)\,=\,\frac{1}{\sqrt{n}}\,\int^n_0\,\phi_j(t)\,\d \xi_t
 \end{equation}
is gaussian with non-degenerate covariance matrix
$B_\zs{k,n}=\E\,\zeta(n)\zeta^\prime(n)$.}\\[2mm]
\noindent
$\C_2)$  {\sl The maximal eigenvalues of  matrices
$B_\zs{k,n}$
satisfy the following inequality
$$
\sup_\zs{k\ge 1}\,
\sup_\zs{n\ge 1}\,
\lambda_\zs{\max}(B_\zs{k,n})\,\le\,\lambda^*\,,
$$
 where
$\lambda^*$ is some known positive constant.}

 Processes \eqref{1.2} and \eqref{1.3} in Examples~\ref{Ex.1.1}--\ref{Ex.1.2}, 
as is shown in Lemmas~\ref{Le.A.0}--\ref{Le.A.0-1},
satisfy condition
$\C_\zs{1})$.
Condition
$\C_\zs{2})$ is satisfied for process \eqref{1.2} with $\lambda^*=2$.
Condition
$\C_\zs{2})$  holds also for process \eqref{1.3} provided that the
value of
vector $\theta$ belongs to the following compact set
 \begin{equation}\label{2.1-1}
K_\zs{\delta}=
\left\{\theta\in\cA\,:\,\max_\zs{1\le i\le q}\,
\mbox{Re}\,\lambda_\zs{i}(\theta)\le -\delta\,,\quad |A(\theta)|\le \delta^{-1}
\right\}\,,
 \end{equation}
where $0<\delta<1$ is a known constant;  $|\cdot|$
stands for
the euclidean norm of  matrix. Under this assumption
process \eqref{1.3} satisfies  condition $\C_\zs{2})$
with
\begin{equation}\label{2.1-2}
\lambda^*=\lambda^*(\delta)=
\frac{2}{\delta^2}
F^*(\delta)\,J^*(\delta)\,,
 \end{equation}
where
$$
F^*(\delta)=
\frac{q}{2\delta}
+
\frac{2q}{\delta^3}
\sum^{q-1}_\zs{j=1}\frac{(2j)!}{(j!)^2\delta^{4j}}
\quad\mbox{and}\quad
J^*(\delta)=
\frac{1}{\delta}+\frac{2}{\delta^2}
\sum^{q-1}_\zs{j=1}\frac{2^{j}}{\delta^{2j}}\,.
$$
Let $\bbn$ be the set of positive integer numbers, i.e.
$\bbn=\{1,2,\ldots\}$. Denote by $\cM$  some finite set of
finite subsets of $\bbn$ and by $(\cD_\zs{m})_\zs{m\in\cM}$ a
family of linear subspaces of $\cX$ such that
$$
 \cD_\zs{m}=\{x\in
\cX\,:\,x=\sum_{j\in m}\lambda_j\phi_j\,,
 \lambda_j\in\bbr\}\,.
$$

 Let $d_m=\dim{\cD_\zs{m}}$ be the number of
elements in a subset $m$. Denote by $S_m$ the projection of $S$ on
$\cD_\zs{m}$, i.e.
\begin{equation}\label{2.2}
S_m=\sum_\zs{j\in m}\alpha_j\phi_j\,,\ \ \  \alpha_j=(S,\phi_j)\,.
\end{equation}

To estimate the function $S$ in \eqref{1.1} we will apply a
general model selection approach. It requires first to choose some
class of projective estimators $\wt{S}_m$ of $S_m$, which may
be any measurable functions of observations $(y_t)_\zs{0\le t\le
n}$ taking on values in $\cD_\zs{m}$. For example, one can take
 the LSE $\wh{S}_m$ of $S$, which is the minimizer, with respect to $x\in \cD_\zs{m},$
of the quantity
\begin{equation}\label{2.3}
\gamma_n(x)=\|x\|^2\,-\,2\,\frac{1}{n}\, \int^n_0\,x(t)\,\d y_t\,
\end{equation}
and has the form
\begin{equation}\label{2.4}
\wh{S}_m=\sum_\zs{j\in m}\,\wh{\alpha}_j\,\phi_j\,,\quad
\wh{\alpha}_j\,=\,\frac{1}{n}\,\int^n_0\,\phi_j(t)\,\d y_t\,.
\end{equation}
 Let $(l_\zs{m})_\zs{m\in \cM}$
be a sequence of prior weights  such that $l_\zs{m}\ge 1$ 
for all $m\in\cM$. We set
\begin{equation}\label{2.5}
l^{*}=\sum_{m\in\cM} e^{-l_\zs{m}\,d_\zs{m}}\,.
\end{equation}

Further one needs a penalty term on the set $\cM$. We take it in the form
suggested in Birg\'e and Massart (2001).
We define the penalty term as
\begin{equation}\label{2.6}
P_\zs{n}(m)\,=\,\rho\, \frac{l_\zs{m}\,d_m}{n}
\quad\mbox{with}\quad
 \rho=\,4\lambda^*\,\frac{z^2_\zs{*}}{z_\zs{*}-1}\,,
\end{equation}
where $z_\zs{*}$ is the maximal root of the equation
$\ln z=z-2$ which is approximately equal to
$z_\zs{*}\approx 3,1462$.

Minimizing the penalized empirical contrast
$\gamma_n(\wt{S}_m)+P_\zs{n}(m)$ with respect to
$m\in\cM$ one finds
\begin{equation}\label{2.7}
\wt{m}=\mbox{argmin}_{m\in\cM}\{\gamma_n(\wt{S}_m)+P_\zs{n}(m)\}
\end{equation}
and obtains the model selection procedure
$\wt{S}_\zs{\wt{m}}$ corresponding to a specific class of
projective estimators $(\wt{S}_\zs{m})_\zs{m\in\cM}$. For
the  LSE family $(\wh{S}_\zs{m})_\zs{m\in\cM}$, this yields
$\wh{S}_\zs{\wh{m}}$ with
\begin{equation}\label{2.8}
\wh{m}=\mbox{argmin}_{m\in\cM}\{\gamma_n(\wh{S}_m)+P_\zs{n}(m)\}\,.
\end{equation}
 Our first result is  the following.
\begin{theorem}\label{Th.2.1}
Assume that the conditions $\C_1)$--$\C_2)$ are fulfilled for the noise in
\eqref{1.1}.
Then for any class of projective estimators $(\wt{S}_\zs{m})_\zs{m\in\cM}$
the general model selection procedure
 $\wt{S}_{\wt{m}}$ satisfies the
following oracle inequality
\begin{equation}\label{2.9}
\E_\zs{S}\,\|\wt{S}_{\wt{m}}-S\|^2\,
\le\,\inf_{m\in\cM}\,\wt{a}_\zs{m}(S)\,
+\frac{\lambda^{*}\tau_\zs{0}}{n}\,,
\end{equation}
where $\E_\zs{S}$ denotes the expectation with respect to the
distribution of \eqref{1.1} given $S$ ,
$$
\wt{a}_\zs{m}(S)=3\E_\zs{S}\|\wt{S}_m-S\|^2
+16\lambda^{*}z_\zs{*}\,\frac{d_\zs{m}\,l_\zs{m}}{n}
\quad\mbox{and}\quad
\tau_\zs{0}=\frac{16 l^{*} z_\zs{*}}{z_\zs{*}-1}
\,.
$$
\end{theorem}
\noindent The proof of Theorem~\ref{Th.2.1} is given in the
Appendix.
\begin{remark}\label{Re.2.1}
It will be noted that
the choice of the coefficient $\rho$ in the penalty term
\eqref{2.6}, as will be shown in the
proof of the theorem,
provides the minimal value of the principal term
$\wt{a}_\zs{m}(S)$.
\end{remark}

Now we will find the upper bound \eqref{2.9} for the LSE model
selection procedure $\wh{S}_\zs{\wh{m}}$ defined by
\eqref{2.4} and \eqref{2.8}. To this end we have to calculate the
accuracy of $\wh{S}_m$ for any $m\in\cM$. We have
\begin{align*}
\E_\zs{S}\,
\|\wh{S}_\zs{m}-S\|^2&=\|S_\zs{m}-S\|^2+\E_\zs{S}\,\|\wh{S}_\zs{m}-S_\zs{m}\|^2\\
&=\|S_\zs{m}-S\|^2+\sum_\zs{j\in m}\,\E_\zs{S}\,(\wh{\alpha}_\zs{j}- \alpha_\zs{j} )^2\,,
\end{align*}
where $S_\zs{m}$ is given in \eqref{2.2}.
 Moreover, the condition $\C_\zs{2})$ yields
$$
\E_\zs{S}\,(\wh{\alpha}_j\,-\alpha_j)^2
=\frac{1}{n^2}\,\E_\zs{S}\,
\left(
\int^n_0\,\phi_j(t)\,\d \xi_t
\right)^2
\le\,\frac{\lambda^{*}}{n}\,.
$$
Therefore
$$
\E_\zs{S}\,
\|\wh{S}_\zs{m}-S\|^2
\,\le\,\|S_\zs{m}-S\|^2+\lambda^{*}\frac{d_\zs{m}}{n}\,.
$$
Thus, we obtain the following result.
\begin{corollary}\label{Co.2.1}
 Under the conditions $\C_\zs{1})$ and 
$\C_\zs{2})$ the LSE model selection procedure
$\wh{S}_{\wh{m}}$, defined by \eqref{2.4} and \eqref{2.8}, satisfies the inequality
\begin{equation}\label{2.10}
\E_\zs{S}\,\|\wh{S}_{\wh{m}} -S\|^2\,
\le\,\inf_{m\in\cM}\,\wh{a}_\zs{m}(S)\,
+\frac{\lambda^{*}\tau_\zs{0}}{n}\,,
\end{equation}
where
$$
\wh{a}_\zs{m}(S)=3\|S_\zs{m}-S\|^2+
\tau_\zs{1}\lambda^{*}\,
\frac{d_\zs{m}l_\zs{m}}{n}\,,
\quad
\tau_1=3+16 z_\zs{*}
\,.
$$
\end{corollary}
Consider the upper bound in \eqref{2.10} in more detail for the model \eqref{1.1}-\eqref{1.2}.
\begin{corollary}\label{Co.2.2} 
For the model \eqref{1.1}--\eqref{1.2} the LSE model selection procedure
$\wh{S}_{\wh{m}}$, defined by \eqref{2.4} and \eqref{2.8} with
$\lambda^*=2$  satisfies, for any $\theta\le 0$, the inequality
\begin{equation}\label{2.11}
\E_\zs{S}\,\|\wh{S}_{\wh{m}}-S\|^2\,
\le\,\inf_{m\in\cM}\,\wh{b}_m(S)\,
+\,\frac{2\tau_0}{n}\,,
\end{equation}
where $\tau_\zs{0}$ is given in \eqref{2.9},
$$
\wh{b}_m(S)=3\|S_m-S\|^2+2
\tau_1\,\frac{d_\zs{m}\,l_\zs{m}}{n}\,.
$$
\end{corollary}
\begin{remark}\label{Re.2.2}
It will be noted that
 for the model \eqref{1.1}--\eqref{1.2} the LSE
 model selection procedure satisfies the oracle inequality uniformly
in the nuisance parameter $\theta$ including the boundary of the
stationarity region of the Ornstein-Uhlenbeck process, i.e.
$\theta=0$.
\end{remark}

\section{The improvement of LSE.}

In this section we consider a special case of the model
\eqref{1.1}--\eqref{1.2} when $\theta=0$, i.e.
\begin{equation}\label{3.0}
\d y_t\,=\,S(t)\,\d t\,+\,\d w_t\,.
\end{equation}
By applying  the improvement method proposed in Fourdrinier and Pergamenshchikov (2007)
 we will show
that the upper bound in the oracle inequality can be
lessened by a proper choice of the projective estimators.
 Let us introduce a class of estimators of the form
\begin{equation}\label{3.1}
S^*_m(t)=\wh{S}_m(t)+\Psi_m(\wh{S}_m)(t)\,.
\end{equation}
Here $\Psi_m$ is a function from $\bbr^d$ into  $\cD_\zs{m}$, i.e.
\begin{equation}\label{3.2}
\Psi_m(x)(t)=\,\sum_\zs{j\in m}\,v_j(x)\,\phi_j(t)\,,
\quad
x\in\bbr^d
\end{equation}
and $(v_j(\cdot))_\zs{j\in m}$ are
$\bbr^d\to\bbr$ functions such that
$\E_\zs{S}\,v^2_j(\wh{\alpha})<\infty$, where\\
 $\wh{\alpha}=(\wh{\alpha}_j)_\zs{j\in m}$ is the vector with the
components $\wh{\alpha}_j$ defined in \eqref{2.4}. The 
functions $v(x)=(v_j(x))_\zs{j\in m}$ will be specified below. Let
\begin{equation}\label{3.3}
\Delta_m(S)\,=\,\E_\zs{S}\,\|S^*_m-S_m\|^2\,-\,\E_\zs{S}\,\|\wh{S}_m-S_m\|^2\,.
\end{equation}
It is easy to check that
\begin{equation}\label{3.4}
\Delta_m(S)\,=\,2\,\E_\zs{S}\,(\Psi_m\,,\,\wh{S}_m-S_m)\,+\,\E_\zs{S}\,\|\Psi_m\|^2\,.
\end{equation}
 This function can be found explicitly  for the model \eqref{3.0}.
\begin{lemma}\label{Le.3.1}
Let $S^*_m$ be defined by \eqref{3.1}--\eqref{3.2} with
continuously differentiable functions $v_j$ such that
$\E_\zs{S}\,v^2_j(\wh{\alpha})<\infty$. Then
$\Delta_m(S)\,=\,\E_\zs{S}\,L(\wh{\alpha})$,
where
\begin{equation}\label{3.5}
L(x)\,=\,\frac{2}{n}\,\mbox{\em div}\,v(x)\,+\,\|v(x)\|^2\,.
\end{equation}
\end{lemma}
\noindent {\bf Proof}. From \eqref{2.4}, \eqref{3.1}, one has
\begin{align*}
\wh{\alpha}_j\,=\,\frac{1}{n}\,\int^n_0\,\phi_j(t)\,\d y_t
\,=\,\alpha_j\,+\,\frac{1}{n}\,\int^n_0\,\phi_j(t)\,\d w_t\,,
\end{align*}
where $\alpha_j=\int^1_0\,\phi_j(t)S(t)\,\d t$. Therefore
the vector $\wh{\alpha}=(\wh{\alpha}_j)_\zs{j\in m}$
has a normal distribution $\cN(\alpha,n^{-1}I_\zs{d})$, where
$\alpha=(\alpha_\zs{j})_\zs{j\in m}$ and $I_d$ is the unit matrix of order
$d$. This enables one to find the explicit expression for the
first term in the right-hand side of \eqref{3.4}.
Indeed,
\begin{align*}
J\,&=\,\E_\zs{S}\,(\Psi_m,\wh{S}_m-S_m)\,
=\,\E_\zs{S}\,\sum_\zs{j\in m}\,v_j(\wh{\alpha})\,(\wh{\alpha}_j-\alpha_j)\\
&=\,\E_\zs{S}\,(v(\wh{\alpha})\,,\,\wh{\alpha}-\alpha)\,
=\,\int_\zs{\bbr^d}\,(v(u)\,,\,u-\alpha)\,g(\|u-\alpha\|^2)\,\d u\,,
\end{align*}
where
\begin{equation}\label{3.6}
g(a)\,=\,\left(\frac{n}{2\pi}\right)^{d/2}\,e^{-n a/2}\,.
\end{equation}
Making the spherical changes of the variables yields
\begin{align*}
J\,&=\,\int^\infty_0\,\int_\zs{\cS_\zs{r,d}}\,(v(u)\,,\,u-\alpha)\,
\nu_\zs{r,d}(\d u)\,g(r^2)\,\d r\\
&=\,\int^\infty_0\,\int_\zs{\cS_\zs{r,d}}\,(v(u)\,,\,e(u))\,
\nu_\zs{r,d}(\d u)\,r\,
g(r^2)\,\d r\,,
\end{align*}
where $\cS_\zs{r,d}=\{u\in\bbr^d\,:\,\|u-\alpha\|=r\}$,
$\nu_\zs{r,d}(\cdot)$ is the superficial measure on the sphere
$\cS_\zs{r,d}$  and $e(u)=(u-\alpha)/\|u-\alpha\|$ is a normal
vector to
 this sphere.
 By applying the Ostrogradsky--Stokes divergence  theorem we obtain that
\begin{align*}
J=\int^\infty_0\,\int_\zs{\cB_\zs{r,d}}\,\mbox{div}\,v(u)\,\d u
\,r\,g(r^2)\,\d r
\end{align*}
with $\cB_\zs{r,d}=\{u\in\bbr^d\,:\,\|u-\alpha\|\le r\}$. By the
Fubini theorem and the definition of $g$ in (\ref{3.6}) one gets
\begin{align*}
J\,&=\,\frac{1}{2}\int_\zs{\bbr^d}\,
\int^\infty_\zs{\|u-\alpha\|^2}\,g(a)\,\d a
\,\mbox{div}\,v(u)\,\d u\\
&=\,\frac{1}{n}\int_\zs{\bbr^d}\,g(\|u-\alpha\|^2)\,
\,\mbox{div}\,v(u)\,\d u\,=\,
\frac{1}{n}\,\E_\zs{S}\,\mbox{div}\,v(\wh{\alpha})\,.
\end{align*}
This leads to the assertion of Lemma~\ref{Le.3.1}.
\endproof

 In particular for $d_\zs{m}>2$, if
one takes 
$$ 
v(u)=-\frac{(d_m-2)u}{n \|u\|^2}\,,
\quad\mbox{then}\quad
 L(u)=-\frac{(d_m-2)^2}{n^2 \|u\|^2}\,,
 $$
and hence, $\Delta_m(S)\,<\,0$, that is, the estimate \eqref{3.1}
outperforms the least squares esimate \eqref{2.4} in the
approximation of $S_m.$ This allows to improve the model selection
procedure by making use of the estimates \eqref{3.1} instead of
the least squares $\{\wh{S}_m\}.$ As a direct consequence of
Theorem ~\ref{Th.2.1}, one obtains the following result for the
improved model selection procedure $S^*_\zs{m^*}$.

\begin{theorem}\label{Th.3.1}
For the model \eqref{3.0}
the improvement model
selection procedure $S^*_{m^*}$ defined by
(\ref{2.7})
with $\wt{S}_\zs{m}=S^*_\zs{m}$ and $\lambda^*=2$
satisfies the inequality
\begin{equation}\label{3.8}
\E_\zs{S}\,\|S^*_{m^*}-S\|^2
\le\,\inf_{m\in\cM}\,u^*_m(S)
+\frac{2\tau_\zs{0}}{n}\,,
\end{equation}
where $\tau_\zs{0}$ is given in \eqref{2.9},
$$
u^*_m(S)\,=\,3\,\E_\zs{S}\,\|S^*_m-S\|^2\,+\,
32\,z_\zs{*}\,\frac{d_\zs{m}\,l_\zs{m}}{n}\,.
$$
\end{theorem}

\medskip

\section{Asymptotic estimation}

\subsection{The risk upper bound}

In this  section we consider the asymptotic estimation problem for
the model \eqref{1.1}.  To this end,
we additionally assume that all functions in the orthonormal system
$(\phi_j)_\zs{j\ge 1}$ are $1$-periodic and the
unknown function $S(\cdot)$ in the model \eqref{1.1} belongs to
the following functional class
\begin{equation}\label{4.1}
\Theta_\zs{\beta,r}\,=\,\{S\in \cC(\bbr)\cap\cX
\,:\,
\max_{n\ge 1}n^{2\beta}\,\varsigma_\zs{n}(S)\,\le r^2\}\,.
\end{equation}
Here  $\cC(\bbr)$ denotes the set of all continuous
$\bbr\to\bbr$ functions and
\begin{equation}\label{4.1-1}
\varsigma_\zs{n}(S)=\sum^\infty_{j=n}\,s^2_j\,,
\end{equation}
where $(s_j)_\zs{j\ge 1}$ are the Fourier coefficients for the  basis
$(\phi_j)_\zs{j\ge 1}$, i.e.
$$
s_j=(S,\phi_j)=\int^1_0\,S(t)\,\phi_j(t)\,\d t\,;
$$
 $\beta>0$ and $r>0$ are  unknown constants.

Similarly to Galtchouk and Pergamenshchikov (2006) we define now the risk for an estimator $\wt{S}_n$
(a measurable function of the observation $(y_t)_\zs{0\le t\le n}$
in \eqref{2.1}) as follows
\begin{equation}\label{4.2}
\cR_n(\wt{S}_n,\beta)\,=\,\sup_{S\in\Theta_\zs{\beta,r}}\,\sup_\zs{Q\in\cP_\zs{\kappa}}\,
\E_\zs{S,Q}\|\omega_n(\wt{S}_n-S)\|^2\,,
\end{equation}
where $\omega_n=\omega_\zs{n}(\beta)=n^{\frac{\beta}{2\beta+1}}$.
Here $\cP_\zs{\kappa}$ is some class of distributions $Q$ (in the space
$\C[0,+\infty)$) of the noise process $(\xi_\zs{t})_\zs{t\ge 0}$
satisfying conditions $\C_\zs{1})$ and $\C_\zs{2})$ with
$\lambda^{*}=\lambda^{*}(Q)\le \kappa<\infty$ for some known fixed parameter
$\kappa$. In addition, this class is assumed to include the Wiener
distribution $Q_\zs{0}$. The second index in $\E_\zs{S,Q}$
denotes that the expectation is taken with respect to
the distribution of the process \eqref{1.1} corresponding
to the noise distribution $Q$.

Note that for the model \eqref{1.1}--\eqref{1.2}
$\cP_\zs{\kappa}$ is the class of distributions of the
processes \eqref{1.2} with $\theta\le 0$. In this case
$\kappa=2$. For the model \eqref{1.1}--\eqref{1.3}
$$
\cP_\zs{\kappa}=\cQ_\zs{\delta}\cup \{Q_\zs{0}\}\,,
$$
where $\cQ_\zs{\delta}$ is the family of distributions
of the processes $\eqref{1.3}$ of order $q\ge 2$ with the
parameters belonging to the set \eqref{2.1-1} for some $0<\delta<1$.
In this case $\kappa=\max(2,\lambda^{*}(\delta))$, where
$\lambda^{*}(\delta))$ is given in \eqref{2.1-2}.

We will apply the ordered variable model selection procedure
(see, Barron et al. (1999), p. 315), for which $\cM=\{m_1\,,\ldots\,,m_n\}$
with $m_i=\{1,\ldots,i\}$, therefore $d_{m_i}=i$. Then
$$
\cD_\zs{m}=\{x\in \cX_n\,:
\,x=\sum^i_{j=1}\,\alpha_j\,\phi_j\,,
\alpha_j\in \bbr\}\,.
$$
 For the ordered variable model selection procedure
 one can take $l_\zs{m}=1$ for all $m\ge 1$  and find
$$
l^{*}=\sum^{n}_{i=1}e^{-d_\zs{m_\zs{i}}}\le \frac{1}{e-1}\,.
$$
In the sequel we denote by $\wh{S}^{\kappa}_{\wh{m}}$ the 
LSE model selection procedure \eqref{2.4}, \eqref{2.8}
 replacing $\lambda^*$ by $\kappa$. Now we will show that the risk \eqref{4.2}
for this procedure
 is finite.
\begin{theorem}\label{Th.4.1} The estimator $\wh{S}^{\kappa}_{\wh{m}}$
satisfies the following asymptotic inequality
\begin{equation}\label{4.3}
\limsup_{n\to\infty}\,\sup_\zs{\beta>0}\,\cR_n(\wh{S}^{\kappa}_{\wh{m}},\beta)<\infty\,.
\end{equation}
\end{theorem}
\proof Taking into account \eqref{2.10} one gets,
for any $Q\in \cP_\zs{\kappa}$,
\begin{align*}
\E_\zs{S,Q}\|\wh{S}^{\kappa}_{\wh{m}}-S\|^2\,&\le\,\inf_{m\in\cM_n}
\left(3\,\|S_\zs{m}-S\|^2\,+\,\tau_1\,\kappa\frac{d_\zs{m}}{n}\right)\,
+\frac{\kappa \tau_\zs{0}}{n}\\
&\le\,\inf_{1\le i\le n}
\left(3\,\|S_\zs{m_i}-S\|^2\,+\,\tau_1\,\kappa\,\frac{i}{n}\right)\,
+\frac{\kappa \tau_0}{n}\,.
\end{align*}
Further, for any function $S$ from
$\Theta_\zs{\beta,r},$ one has
\begin{align*}
\|S_\zs{m_i}-S\|^2=\sum^\infty_{j=i+1}\,s^2_j\,
\le\,r^2\,i^{-2\beta}\,.
\end{align*}
Therefore  for each $1\le i\le n$
\begin{align*}
\sup_\zs{Q\in\cP_\zs{\kappa}}\,
\E_\zs{S,Q}\|\wh{S}^{\kappa}_{\wh{m}}-S\|^2\,
\le\,
\left(3\,r^2\,i^{-2\beta}\, +\,\tau_1\kappa \frac{i}{n}
\right)+\frac{\tau_0\kappa}{n}\,.
\end{align*}
Substituting
$i=i_n=[n^{\frac{1}{2\beta+1}}]+1$
leads to \eqref{4.3}.
\endproof

\medskip

\subsection{The risk lower bound}

Now we study the lower bound for the risk \eqref{4.2}.
 We assume that  the orthogonal system $(\phi_\zs{j})_\zs{j\ge 1}$
in \eqref{4.1} is
trigonometric, i.e.
\begin{equation}\label{5.1}
\phi_1(x)\,\equiv\, 1\,,
\quad\mbox{and for}\quad j\ge 2\quad
\phi_j(x)\,=\,\sqrt{2}\,\mbox{Tr}_j(2\pi[j/2] x)\,,
\end{equation}
where $\mbox{Tr}_\zs{j}(x)=\cos x$ for even $j\ge 1$
and $\mbox{Tr}_\zs{j}(x)=\sin x$ for odd $j\ge 1$; $[a]$ is the integer part of $a$.
\medskip
\begin{theorem}\label{Th.5.1}
 The lower bound of the risk \eqref{4.2} over all
estimates is strictly positive, i.e. for any $\beta\ge 1$
\begin{equation}\label{5.2}
\underline{\lim}_{n\to\infty}\inf_{\wt{S}_n}\cR_n(\wt{S}_n,\beta)>0\,.
\end{equation}
\end{theorem}
\noindent {\em Proof}. In order to show \eqref{5.2}
it suffices to check this inequality for the model \eqref{3.0}, i.e.
that
for any $\beta\ge 1$
\begin{equation}\label{5.2-1}
\underline{\lim}_{n\to\infty}\,\inf_{\wt{S}_n}\,\sup_{S\in\Theta_\zs{\beta,r}}\,
\E_\zs{S,Q_\zs{0}}\|\omega_n(\wt{S}_n-S)\|^2
\,>\,0\,.
\end{equation}
 To this end we apply the method proposed in
 Fourdrinier and Pergamenshchikov (2007) to our case.
First, we construct an auxiliary parametric class
of functions in the set $\Theta_\zs{\beta,r}$. Let
$\beta=k+\alpha$ with $k=[\beta]$ and $0\le \alpha<1$.
Let $V(\cdot)$ be  $k+1$ times continuously differentiable function 
such that $V(u)=0$ for $|u|\ge 1$ and
$\int^{1}_\zs{-1}\,V^2(u)\,\d u=1$. Let
$m=[n^{\frac{1}{2\beta+1}}]$  and
  $\Gamma_\zs{\delta}$ be a cube in $\bbr^m$ of the form
$$
\Gamma_\zs{\delta}\,=\,
\{z=(z_1,\ldots,z_m)'\in\bbr^m\,:\, |z_i|\le \delta\,,\quad 1\le j\le m \}\,,
$$
where $\delta=\nu/\omega_\zs{n},$ $\nu>0$. Now, viewing the
function $V(\cdot)$ as a kernel, one introduces
a parametric class of $1$-periodic functions
$(S_z)_\zs{z\in\Gamma_\zs{\delta}}$
where
\begin{equation}\label{5.3}
S_z(t)=\sum^m_{j=1}\,\,z_{j}\psi_\zs{j}(t)\,,\quad 0\le t\le 1;
\end{equation}
$(\psi_j(t))_\zs{1\le j\le m}$ are $1$-periodic functions
 defined on the interval $[0,1]$ as\\ $\psi_j(t)=\,V\left(\frac{t-a_j}{h}\right)$ with
$h=1/2m$ and $a_j=(2j-1)/2m$.

It will be observed that, for $0\le i\le k-1$ and $z\in\Gamma_\delta,$
$$
\sup_\zs{0\le t\le 1}\,|S^{(i)}_z(t)|
\le\, 2^i\,\sup_\zs{|a|\le 1}\,|V^{(i)}(a)|\,\frac{\nu}{n^{(\beta-i)/(2\beta+1)}}
\to 0
\quad\mbox{as}\quad
n\to\infty\,.
$$
In order to check
 the second condition in \eqref{A.9}, we estimate the increment of $k$th
 derivative of $S_\zs{z}(\cdot)$.
For any $0\le s\,,\, t\le 1$ and  $z\in\Gamma_\delta$, one  has
\begin{equation}\label{5.4}
|S^{(k)}_z(t)\,-\,S^{(k)}_z(s)|
\le
\frac{\nu}{\omega_n\,h^k}\,\Delta_\zs{m}\,,
\end{equation}
where
$$
\Delta_\zs{m}\,=\,\sum^m_\zs{j=1}\,
\left|V^{(k)}\left(\frac{t-a_j}{h}\right)\,-\,V^{(k)}\left(\frac{s-a_j}{h}\right)\right|\,.
$$
If $s$ and $t$ belong to the same interval, that is,
 $a_\zs{j_0}-h\le s\,\le\,t\le a_\zs{j_0}+h$,
 then putting $V^*=\sup_\zs{|a|\le 1}\,|V^{(k+1)}(a)|$
 one obtains
\begin{align}\nonumber
\Delta_\zs{m}\,&=\,
\left|V^{(k)}\left(\frac{t-a_\zs{j_0}}{h}\right)\,
-\,V^{(k)}\left(\frac{s-a_\zs{j_0}}{h}\right)\right|\\ \label{5.5}
&\le\,2\,V^*\,\frac{|t-s|}{2h}\,
\le\,2^{1-\alpha}\,V^*\,\frac{|t-s|^\alpha}{h^{\alpha}}
\end{align}
for each $0\le \alpha\le 1$. If $s$ and $t$ belong to different
intervals, that is,
$$
a_\zs{j_0}-h\le s\le a_\zs{j_0}+h\,\le\,a_\zs{j_1}-h\le t\le
a_\zs{j_1}+h\,,\quad
j_0<j_1\,,
$$
then setting $s^*=a_\zs{j_0}+h$ and $t_*=a_\zs{j_1}-h,$ similarly to
(\ref{5.5}), one gets
\begin{align*}
\Delta_\zs{m}\,&=\,
\left|V^{(k)}\left(\frac{t-a_\zs{j_1}}{h}\right)\,-
\,V^{(k)}\left(\frac{t_*-a_\zs{j_1}}{h}\right)\right|\\
&+
\left|V^{(k)}\left(\frac{s-a_\zs{j_0}}{h}\right)
-V^{(k)}\left(\frac{s^*-a_\zs{j_0}}{h}\right)\right|\\
&\le\,2^{1-\alpha}\,V^*\,\frac{1}{h^{\alpha}}\,
\left(|t-t_*|^{\alpha}\,+\,|s-s^*|^{\alpha}\right)\,
\le\,2^{2-\alpha}\,\frac{1}{h^\alpha}\,V^*\,|t-s|^\alpha\,.
\end{align*}
From here and \eqref{5.4}--\eqref{5.5} we come to the estimate
$$
|S^{(k)}_z(t)\,-\,S^{(k)}_z(s)|\,\le \,2^{2+k}\,\nu\,V^*\,
\frac{m^{\beta}}{\omega_n}\,|t-s|^\alpha\,.
$$

Therefore (see Lemma~\ref{Le.A.2}  in Appendix \ref{Se.A.4}) there exist
$\nu>0$ and $n_\zs{0}\ge 1$ such that $S_z\in\Theta_{\beta,r}$ for all
$z\in\Gamma_\delta$ and $n\ge n_\zs{0}$. Further, we introduce
the prior distribution on $\Gamma_\zs{\delta}$ with the density
$$
\pi(z)=
\pi(z_1,\ldots,z_m)=\prod_{l=1}^m\,\pi_\zs{l}(z_l)\,,
\quad \pi_\zs{l}(u)=\frac{1}{\delta}
G\left(\frac{u}{\delta}\right)\,.
$$
The function
$G(u)=G_\zs{*}e^{-\frac{1}{1-u^2}}$ for $|u|\le 1$ and
$G(u)=0$ for $|u|\ge 1$,
where $G_\zs{*}$ is a positive constant such that $\int^1_{-1}G(u)\,\d u=1$.


Let $\wt{S}_n(\cdot)$ be an estimate of $S(\cdot)$ based on observations
$(y_t)_\zs{0\le t\le n}$ in \eqref{1.1}. Then for any $n\ge n_\zs{0}$  we 
can estimate with below the supremum in \eqref{5.2-1} as
\begin{align*}
\sup_{S\in\Theta_{\beta,r}}\E_\zs{S,Q_\zs{0}}
\|\wt{S}_n-S\|^2
>
\int_{\Gamma_\zs{\delta}}\,
\E_\zs{S_\zs{z},Q_\zs{0}}\|\wt{S}_n-S_z\|^2\,
\pi(z)\,\d z\,.
\end{align*}
Moreover, by the definition of $S_\zs{z}$, we obtain
\begin{align*}
\|\wt{S}_n-S_z\|^2
&\ge
\sum^m_{l=1}(\wt{z}_l-z_l)^2\int^1_0\psi^2_l(t)\,\d t
=h \sum^m_{l=1}(\wt{z}_l-z_l)^2\,,
\end{align*}
where $\wt{z}_l=\,\int^1_0\,\wt{S}_n(x)\,\psi_l(x)\,\d x/\|\psi_l\|^2$.
Therefore,
\begin{equation}\label{5.6}
\sup_{S\in\Theta_\zs{\beta,r}}\E_\zs{S,Q_\zs{0}}\,\|\wt{S}_n-S\|^2
\ge
h\,\sum^m_{l=1}\Lambda_l\,,
\end{equation}
where $\Lambda_l=\int_{\Gamma_\delta} \E_\zs{S_\zs{z}}(\wt{z}_l-z_l)^2 \pi(z)\d z$.
To apply now lemma~\ref{Le.A.2-1} we note that in this case
$$
\zeta_l(z)=\int^n_0\,\psi_l(t)\,\d y_t\,
-\,\int^n_0\,S_z(t)\,\psi_l(t)\d t
$$
and, therefore, $A_\zs{l}=
\E_\zs{S_\zs{z},Q_\zs{0}}\,\zeta^2_l(z)=\int^n_0\,\psi^2_l(t)\,\d t\,
=\,nh$.
Moreover, in this case
$$
B_l=
\int^\delta_{-\delta}
\frac{(\dot{\pi}_\zs{l}(u))^2}{\pi_\zs{l}(u)}\d u
=\delta^{-2} I_\zs{G}
\quad ;\quad
I_\zs{G}=8\int^1_0u^2(1-u^2)^{-4}G(u)\,\d u\,.
$$
Thus, by the inequality \eqref{A.10-3}, one obtains that
\begin{align*}
\sup_{S\in\Theta_{\beta,r}}\,
\E_\zs{S,Q_\zs{0}}
\|\wt{S}_n-S\|^2
&\ge \frac{1}{2m}\sum^m_{l=1}
\frac{1}{nh+\omega^2_n \nu^{-2}I_\zs{G}}\\
&
=\,\frac{1}{2nh+2\omega^2_n \nu^{-2}I_\zs{G}}\,.
\end{align*}
This  immediately implies \eqref{5.2-1}.
Hence Theorem~\ref{Th.5.1}.
\endproof

\section{Estimation  based on discrete data.}

The model selection procedure developed in Section 2 is intended for
continuous time observations. However, in a number of applied problems high
frequency sampling can not be provided.
In this section, we consider the estimation problem for model \eqref{1.1}
on the basis of  observations $(y_\zs{t_j})_\zs{0\le j\le np}$ of the process
$(y_\zs{t})_\zs{t\ge 0}$ at discrete times $t_\zs{j}=j/p$, where $p$ is a
given odd
number. To solve this problem, we will modify the model selection procedure of Section 2.
 Let
$\cX_\zs{p}$ be the set of all $1$-periodic  functions
$x\,:\,\bbr\,\to\,\bbr$ with the scalar product
\begin{equation}\label{6.1}
(x,z)_\zs{p}\,=\,\frac{1}{p}\,\sum^{p}_\zs{j=1}\,x(t_j)\,z(t_j)\,,
\quad x,z\in\cX_\zs{p}\,.
\end{equation}

Let
$(\phi_j)_\zs{1\le j\le p}$ be an orthonormal basis in $\cX_\zs{p}$,
i.e. $(\phi_i,\phi_j)_\zs{p}=0$, if $i\ne j$ and
$\|\phi_i\|^2_\zs{p}=1$. One can use, for example, the
 trigonometric basis \eqref{5.1}.

Assume that the noise $(\xi_t)_\zs{t\ge 0}$ in \eqref{1.1} is such that

\noindent $\C^{*}_1)$ {\sl  The vector
$\zeta^{*}(n)=(\zeta^{*}_\zs{1}(n),\ldots,\zeta^{*}_\zs{p}(n))^\prime$ with
components
 \begin{equation}\label{6.4}
  \zeta^{*}_\zs{l}(n)\,=\,\frac{1}{\sqrt{n}}\,\sum^n_{j=1}\,\phi_\zs{l}(t_j)\,\Delta \xi_\zs{t_j}\,, \ \
\Delta \xi_\zs{t_j}=\xi_\zs{t_j}-\xi_\zs{t_{j-1}}\,,
 \end{equation}
is gaussian with non-degenerate covariance matrix
$B^{*}_\zs{n,p}=\E\,\zeta^{*}(n)(\zeta^{*}(n))^\prime$;}\\[2mm]
\noindent $\C^{*}_2)$ {\sl The maximal eigenvalues of matrices $B^{*}_\zs{n,p}$ are uniformly
bounded :
$$
\sup_\zs{n\ge 1}
\sup_\zs{p\ge 1}\,\lambda_\zs{\max}(B^{*}_\zs{n,p})\,\le\,\lambda^*\,,
$$
where $\lambda^*$ is some known positive constant.
}

Conditions $\C^{*}_\zs{1})$, $\C^{*}_\zs{2})$ are satisfied for
processes \eqref{1.2}, \eqref{1.3} (cf.
Lemmas~\ref{Le.A.0}--\ref{A.0-1}).

 Now we denote by $\cM_\zs{p}$  some set of subsets of $\{1,\ldots,p\}$
and by $(\cD_\zs{m,p})_\zs{m\in\cM_\zs{p}}$  a family of linear
subspaces of $\cX_\zs{p}$ such that
$$
\cD_\zs{m,p}\,=\,\{x\in \cX_\zs{p}\,:\,x=\sum_{j\in m}\,\lambda_j\phi_j\,,
 \lambda_j\in\bbr\}\,.
$$
Let $S_\zs{m,p}$ denote the projection of $S$ on
$\cD_\zs{m,p}$ in $\cX_\zs{p}$
and $\wt{S}_\zs{m,p}$ denote an estimator
of $S_\zs{m,p}$, i.e. a measurable function of the observations
$(y_\zs{t_j})_\zs{0\le j\le np}$ taking on values in $\cD_\zs{m,p}$.
One can use, for example, the LSE $\wh{S}_\zs{m,p}$ for $S_\zs{m,p}$, which
 is defined as
the minimizer with respect to $x\in \cD_\zs{m,p}$ of the distance
$$
\frac{1}{np}\,\sum^{np}_\zs{k=1}\,
\left(\frac{\Delta y_\zs{t_k}}{\Delta t_k}\,-\,x(t_k)\right)^2
$$
that is , the quantity
\begin{equation}\label{6.5}
\gamma_\zs{n,p}(x)=\|x\|^2_p\,-\,2\,\frac{1}{n}\,
\sum^{np}_\zs{k=1}\,x(t_k)\,\Delta y_\zs{t_k}\,
\end{equation}
and has the  form
\begin{equation}\label{6.6}
\wh{S}_\zs{m,p}=\sum_\zs{j\in
m}\,\wh{\alpha}_\zs{j,p}\,\phi_j\,,\quad
\wh{\alpha}_\zs{j,p}\,=\,\frac{1}{n}\,\sum^{np}_\zs{k=1}\,\phi_j(t_k)\,\Delta y_\zs{t_k}\,.
\end{equation}
Let the penalty term $P_\zs{n}(m)$ be defined, as before, by \eqref{2.6}.  Then the model
selection procedure,
corresponding to a family of projective
estimators $(\wt{S}_\zs{m,p})_\zs{m\in\cM_\zs{p}}$,
is defined as $\wt{S}_\zs{\wt{m}_\zs{p},p}$ where
\begin{equation}\label{6.7}
\wt{m}_\zs{p}=\mbox{argmin}_{m\in\cM_\zs{p}}\{\gamma_\zs{n,p}(\wt{S}_\zs{m,p})+P_\zs{n}(m)\}\,.
\end{equation}
 In the case of the LSE family
 $(\wh{S}_\zs{m,p})_\zs{m\in\cM_\zs{p}}$,
it will be $\wh{S}_\zs{\wh{m}_\zs{p},p}$.

As a measure of  accuracy
of the approximation of a
 $1$-periodic function $S$ of continuous argument $t$ by its values on the
$(t_j)_\zs{1\le j\le p}$, we will use the function
\begin{equation}\label{6.9}
H_p(S)\,=\,\frac{1}{p}\,\sum^p_\zs{l=1}\,h^2_\zs{l}(S)\,,\quad
h_\zs{l}(S)\,=\,\frac{1}{\Delta t_l}\int^{t_l}_\zs{t_{l-1}}\,(S(t)-S(t_l))\,\d t\,.
\end{equation}
\medskip
The following theorem gives the oracle inequality for a general
 model selection procedure  $\wt{S}_\zs{\wt{m}_\zs{p},p}$ based on the
 discrete time
observations.
\begin{theorem}\label{Th.6.1}
Assume that the conditions $\C^{*}_1)$--$\C^{*}_2)$ hold.
Then the estimator $\wt{S}_\zs{\wt{m}_\zs{p},p}$
 satisfies the  oracle inequality
\begin{equation}\label{6.10}
\E_\zs{S}\,\|\wt{S}_\zs{\wt{m}_\zs{p},p}-S\|^2_p\le
\inf_{m\in\cM_\zs{p}}\,\wt{a}_\zs{m,p}(S)+8H_p(S)
+\frac{2 \tau_\zs{0}\lambda^{*}}{n}\,,
\end{equation}
where
$$
\wt{a}_\zs{m,p}(S)\,=\,7\E_\zs{S}\,\|\wt{S}_\zs{m,p}-S\|^2_p\,
+\,32\lambda^{*}z_\zs{*}\,\frac{l_\zs{m} d_m}{n}\,.
$$
\end{theorem}

Now we obtain the oracle inequality \eqref{6.10}
for the least square model selection procedure
$\wh{S}_\zs{\wh{m}_\zs{p},p}$.
 To this end, we have to calculate
the estimation accuracy of $\wh{S}_\zs{m,p}$ for
$S_\zs{m,p},$ which is
the projection of $S$
on $\cD_\zs{m,p}$, i.e.
$$
S_\zs{m,p}=\sum_\zs{j\in m}\,\alpha_\zs{j,p}\,\phi_j\,,\quad
\alpha_\zs{j,p}=(S,\phi_j)_p=\frac{1}{p}\,\sum^p_\zs{k=1}\,S(t_k)\,\phi_j(t_k)\,.
$$
First of all, we note that in this case
\begin{align*}
\wh{\alpha}_\zs{j,p}\,-\,\alpha_\zs{j,p}\,=\,
\frac{1}{p}\,\sum^p_\zs{k=1}\,\phi_j(t_k)\,h_k(S)\,+\,\frac{1}{n}\int^n_\zs{0}
\phi_\zs{j,p}(t)\,\d \xi_t\,,
\end{align*}
where
$\phi_\zs{j,p}\,=\,\sum^{np}_\zs{k=1}\,\phi_j(t_k)\,\Chi_\zs{(t_\zs{k-1},t_k]}(t)$.
In view of the condition $\C^{*}_\zs{2})$,
this implies that 
\begin{align*}
\E_\zs{S}\,(\wh{\alpha}_\zs{j,p}\,-\alpha_\zs{j,p})^2&=\,
\frac{1}{p^2}\,\left(\sum^p_\zs{k=1}\,\phi_j(t_k)\,h_k(S)\right)^2
\\
&+
\frac{1}{n^2}\,\E_\zs{S}\,\left(\int^n_0\,\phi_\zs{j,p}(t)\,\d \xi_t\right)^2
\le\,H_p(S)+\frac{\lambda^{*}}{n}\,.
\end{align*}

\begin{corollary}\label{Co.6.1} Under the conditions
$\C^{*}_1)$--$\C^{*}_2)$ the LSE
procedure $\wh{S}_{\wh{m}_\zs{p},p}$
satisfies the inequality
\begin{equation}\label{6.13}
\E_\zs{S}\,\|\wh{S}_\zs{\wh{m}_\zs{p},p}-S\|^2_p
\le\inf_{m\in\cM_\zs{p}}\wh{b}_\zs{m,p}(S)+8H_p(S)
+\frac{2\tau_\zs{0}\lambda^{*}}{ n}\,,
\end{equation}
where
$$
\wh{b}_\zs{m,p}(S)=7\|S_\zs{m,p}-S\|^2_p
+7d_\zs{m}H_\zs{p}(S)
+\lambda^{*}
(7+32 z_*\,l_\zs{m})\frac{d_\zs{m}}{n}\,.
$$
\end{corollary}
Now we consider the estimation problem for the model
\eqref{1.1}--\eqref{1.2} on the basis of discrete data in
 the asymptotic setting.
First, for any $\beta\ge 1$, we set
\begin{equation}\label{6.14}
\cR_\zs{n,p}(\wt{S}_n,\beta)=
\sup_{S\in\Theta_\zs{\beta,r}}\,\sup_\zs{Q\in \cP_\zs{\kappa}}\,
\E_\zs{S,Q}\|\omega_n(\wt{S}_n-S)\|^2_p\,,
\end{equation}
where  the set $\Theta_\zs{\beta,r}$ is defined by \eqref{4.1}
with the use of the trigonometric basis
\eqref{5.1}, $\omega_n=\omega_n(\beta)=n^{\frac{\beta}{2\beta+1}}$
and the set $\cP_\zs{\kappa}$ is defined in \eqref{4.2}.
 As in Section 4, in order to minimize this risk, we apply the least square model selection
procedure $\wh{S}_{\wh{m}_\zs{p},p}$, constructed on
the basis of the trigonometric system \eqref{5.1} with the ordered selection, that is,
$\cM_\zs{p}=\{m_1\,,\ldots,\,m_\zs{p}\}$ with $m_j=\{1,\ldots,j\}$.
 In this case $d_\zs{m_j}=j$ and $l_\zs{m_j}=1$ for $1\le j\le p$.

It is shown in Appendix~\ref{Se.A.6}, that if $p\ge n^{1/2}$, then for any $\varepsilon>0$
\begin{equation}\label{6.15}
\overline{\lim}_\zs{n\to\infty}\,\sup_\zs{\beta\ge 1+\varepsilon}\,
\cR_\zs{n,p}(\wh{S}_\zs{\wh{m}_\zs{p},p},\beta)\,<\,\infty
\end{equation}
and if $p\ge n^{1/2}$, then
for any $\beta\ge 1$
\begin{equation}\label{6.16}
\underline{\lim}_\zs{n\to\infty}\inf_\zs{\wt{S}_n}\,\cR_\zs{n,p}(\wt{S}_\zs{n},\beta)\,>0\,.
\end{equation}
It means that the adaptive estimator $\wh{S}_\zs{\wh{m}_\zs{p},p}$ with the
$p\ge n^{1/2}$ ( in particular, one can take $p=2[n^{1/2}]+1$ )
is optimal in the sense of the risk \eqref{6.14}.

\medskip

\section{Appendix}

\subsection{ Properties of processes \eqref{1.2}--\eqref{1.3}}\label{Se.A.0}

We start with the result for process \eqref{1.2}
which shows that both conditions
$\C_\zs{1})$--$\C_\zs{2})$ and
$\C^{*}_\zs{1})$--$\C^{*}_\zs{2})$ are satisfied.

\begin{lemma}\label{Le.A.0}
Let $(\xi_\zs{t})$ be defined by \eqref{1.2} with $\theta\le 0$,
 $f=(f_\zs{j})_\zs{j\ge 1}$ be a family of  linearly
independent cadlag $1$-periodic functions on $\bbr$,
$I_\zs{t}(f)=\int^t_\zs{0}f(u)\d \xi_\zs{u}$.
Then the matrix
\begin{equation}\label{A.0}
V_\zs{k,n}(f)=(v_\zs{i,j}(f))_\zs{1\le i,j\le k}
\end{equation}
 with elements
$v_\zs{i,j}(f)=\E\, I_\zs{n}(f_\zs{i}) I_\zs{n}(f_\zs{j})$
is positive definite for each $k\ge 1$, $n\ge 1$ and $\theta\le 0$.
Moreover, if $(f_\zs{j})_\zs{j\ge 1}$ is orthonormal, then for any $\theta\le 0$
\begin{equation}\label{A.0-1}
\sup_\zs{k\ge 1}\,\sup_\zs{n\ge 1}\,
\sup_\zs{|z|=1}\,
\frac{1}{n}\,z'V_\zs{k,n}(f)z\le 2\,.
\end{equation}
\end{lemma}
\proof
Assume that for some $n\ge 1$, $k\ge 1$ and $z=(z_\zs{1},\ldots,z_\zs{k})'\in\bbr^k$
$z'V_\zs{k,n}(f)z=0$. Since
$$
z'V_\zs{k,n}(f)z=\,\E\, I^2_\zs{n}(g)\,,
$$
where $g(t)=\sum^k_\zs{j=1}z_\zs{j}f_\zs{j}(t)$,
one gets $I_\zs{n}(g)=\int^t_\zs{0} g(t)\d \xi_t=0$ a.s.
Taking into account that
the distribution of $I_\zs{n}(g)$
for model \eqref{1.2}
 is equivalent to
that of the random variable  $\int^n_\zs{0} g(t)\d w_\zs{t}$
this implies that  $g(t)=0$
for all $t\in [0,n]$. Thus
$z_\zs{1}=\ldots=z_\zs{k}=0$ and to we come the first assertion.
Let us check \eqref{A.0-1}. By applying Ito's formula one obtains
$$
\E\, I^2_\zs{n}(g)=2 \theta \int^n_\zs{0}\,g(t)\,
\E I_\zs{t}(g)\,\xi_\zs{t}\d
t+\int^n_\zs{0} g^2(t)\d t\,,
$$
where
$$
\E I_\zs{t}(g)\,\xi_\zs{t}
=\frac{1}{2}\int^t_\zs{0} g(u) e^{\theta (t-u)}\d u\,.
$$
Therefore,
\begin{equation}\label{A.0-2}
\E I^2_\zs{n}(g)=
\theta \int^n_\zs{0} e^{\theta v}\int^{n}_\zs{v} g(t) g(t-v) \d t\, \d v+
\int^n_\zs{0} g^2(t)\d t\,.
\end{equation}
From here, it follows that
for any $\theta\le 0$
\begin{align*}
z'V_\zs{k,n}(f)z&\le
\int^n_\zs{0} g^2(t)\d t
\left(
1+\theta \int^\infty_\zs{0}\,e^{\theta v}\d v
\right)\\
&\le 2 n\int^1_\zs{0} g^2(t)\d t=2n
\sum^k_\zs{j=1}z^2_\zs{j}=2n\,.
\end{align*}
This completes the proof of
 Lemma~\ref{Le.A.0}.
\endproof
\medskip

\begin{lemma}\label{Le.A.0-1}
Let $(\xi_\zs{t})$ be defined by \eqref{1.3} with $\theta\in\cA$,
 $f=(f_\zs{j})_\zs{j\ge 1}$ be a family of  linearly
independent cadlag $1$-periodic functions on $\bbr$.
Then the matrix \eqref{A.0}
is positive definite for each $k\ge 1$, $n\ge 1$ and $\theta\in\cA$.
Moreover, if $(f_\zs{j})_\zs{1\le j\le k}$ is orthonormal,
then for any $0<\delta<1$ and $\theta\in  K_\zs{\delta}$
\begin{equation}\label{A.0-3}
\sup_\zs{k\ge 1}\,\sup_\zs{n\ge 1}\,
\sup_\zs{|z|=1}\,
\frac{1}{n}\,
z'V_\zs{k,n}(f)z
\le
\lambda^*(\delta)\,,
\end{equation}
where $\lambda^*(\delta)$ is defined in \eqref{2.1-2}.
\end{lemma}
\proof
Let $\eta_\zs{t}$ be process
\eqref{1.3} with zero initial values, i.e.
$$
\d \eta^{(q-1)}_\zs{t}=\left(\sum^q_\zs{j=1}\theta_\zs{j}\eta^{(q-j)}_\zs{t}\right)\d t+\d
w_\zs{t}
$$
and $\eta_\zs{0}=\ldots=\eta^{(q-1)}_\zs{0}=0$. Then $\xi_\zs{t}$
can be written as
\begin{equation}\label{A.0-4}
\xi_\zs{t}=<e^{At}Y>_\zs{q}+\eta_\zs{t}\,,
\end{equation}
where
$<X>_\zs{i}$ denotes the $i$th component of a vector $X$;
$A$ is the matrix defined in \eqref{1.5} and $Y$ is a
gaussian vector in $\bbr^q$ independent of $(\eta_\zs{t})_\zs{t\ge 0}$
with  zero mean and covariance matrix
\begin{equation}\label{A.0-5}
F=\int^\infty_\zs{0}\,e^{Au}\,D_\zs{q}\,e^{A'u}\,\d u\,,
\end{equation}
where $D_\zs{q}$ is $q\times q$ matrix in which 
all elements exept of the $(1,1)$ element are equal to zero
ant the $(1,1)$ element is equal to $1$.
In view of \eqref{A.0-4}, one has
\begin{equation}\label{A.0-5-1}
I_\zs{n}(g)=\zeta+\int^n_\zs{0}g(t)\d \eta_\zs{t}\,,
\end{equation}
where $\zeta=<\int^n_\zs{0}g(t)Ae^{At}\d t Y>_\zs{q}$. Integration by
parts yields
$$
\int^n_\zs{0}g(t)\d \eta_\zs{t}=
\int^n_\zs{0}G_\zs{q-1}(t)\d \eta^{(q-1)}_\zs{t}\,,
$$
where
 $G_\zs{0}(t)=g(t)$ and
$G_\zs{j}(t)=\int^n_\zs{t} G_\zs{j-1}(u)\d u$ for $1\le j\le q-1$.

Now assume that for some $n\ge 1$, $k\ge 1$ and
$z=(z_\zs{1},\ldots,z_\zs{k})'\in\bbr^k$
$z'V_\zs{k,n}(f)z=0$. Since
$$
z'V_\zs{k,n}(f)z=\E\, I^2_\zs{n}(g)=\E\,\zeta^2+
\E\,\left(
\int^n_\zs{0}g(t)\d \eta_\zs{t}
\right)^2\,,
$$
this implies that
$$
\E\,\left(
\int^n_\zs{0}g(t)\d \eta_\zs{t}
\right)^2=
\E\,\left(
\int^n_\zs{0}G_\zs{q-1}(t)\d \eta^{(q-1)}_\zs{t}
\right)^2=0\,.
$$
Taking into account that
 the distribution of the
process $(\eta^{(q-1)}_\zs{t})_\zs{0\le t\le n}$ in $\C[0,n]$
is equivalent to Wiener measure we have
$$
\int^n_\zs{0}G_\zs{q-1}(t)\d w_\zs{t}=0
\quad\mbox{a.s.}
$$
and therefore $G_\zs{q-1}(t)=0$ for all $0\le t\le n$ and hence
$g(\cdot)=0$ and we obtain
$z_\zs{1}=\ldots=z_\zs{k}=0$.
This leads to the first assertion.

Let us show \eqref{A.0-3}. By direct
calculations we find
$$
\E\, I^2_\zs{n}(g)=
2\int^n_\zs{0}<Ae^{Au}F A'>_\zs{q,q}
\left(\int^{n-u}_\zs{0}\,g(u+s)\,g(s)\d s\right)\,
\d u\,.
$$
where $<A>_\zs{i,j}$ denotes the $(i,j)$-th element of  matrix $A$.
By applying the Bunyakovskii-Caushy-Schwartz inequality
one gets
$$
\E I^2_\zs{n}(g)\le 2
\int^{n}_\zs{0}\,g^2(s)\d s
\int^n_\zs{0}|<Ae^{Au}F A'>_\zs{q,q}|
\d u\,.
$$
Since
$$
\int^{n}_\zs{0}\,g^2(s)\d s=n\int^{1}_\zs{0}\,g^2(s)\d s
=n\sum^k_\zs{j=1}z^2_\zs{j}\int^1_\zs{0}f^2_\zs{j}(s)\d s=n\,,
$$
we obtain the estimate
\begin{align*}
\frac{1}{n}\,z'V_\zs{k,n}(f)z\,\le
2\int^\infty_\zs{0}|<Ae^{Au}F A'>_\zs{q,q}|
\d u\,\le 2 |A|^2|F|J(A)
\,,
\end{align*}
where $J(A)=\int^\infty_\zs{0}|e^{Au}|\d u$.
In order to come to \eqref{A.0-3} it remains to use the following
inequality
for  matrix exponents of order $q\ge 2$ (see, for example, in  Kabanov and Pergamenshchikov (2003) on p. 228)
$$
|e^{tB}|\le e^{t\Lambda}
\left(
1+2|B|\sum^{q-1}_\zs{j=1}\frac{1}{j!}(2t|B|)^j
\right)\,,
$$
where $\Lambda=\max_\zs{1\le j\le q}\,\mbox{Re}\, \lambda_\zs{j}$, $\lambda_\zs{j}$ are
eigenvalues of the matrix $B$.

Indeed, from \eqref{A.0-5} for any $A\in K_\zs{\delta}$ we find
that
$$
|F|\le F^*(\delta)
\quad\mbox{and}\quad
J(A)\le J^*(\delta)\,,
$$
where the functions $F^*(\delta)$ and $J^*(\delta)$ are defined in
\eqref{2.1-2}.
 Hence Lemma~\ref{Le.A.0-1}.
\endproof

\subsection{Mean forecast inequality}\label{Se.A.2}

\begin{lemma}\label{Le.A.1} (Galtchouk and Pergamenshchikov (2005))

Let $\alpha$ and $\xi$
 be two positive random variables. Let $\beta$ be a positive real number and
$\{\Gamma_\zs{x}, x\ge 0\}$ be a
family of events such that, for any $x$,
$$
\P(\xi>\alpha+\beta x\,,\,\Gamma_\zs{x})=0.
$$
Assume also that there exists some positive integrable on $\bbr_+$ function
$M(x)$ dominating
$\P(\Gamma_x^c)$.
Then
$$
\E\,\xi\,\le\, \E\,\alpha\, +\, \beta M^*\,,
$$
where $M^*\,=\,\int_0^{\infty}\,M(x)\d x$.
\end{lemma}
\proof
We set $\eta=(\xi-\alpha)_\zs{+}$. Thus
$\E\xi \le \E\alpha + \E\eta$. Moreover,
\begin{align*}
\E\eta&=\int^\infty_\zs{0}\P(\eta>z)\d z=
\int^\infty_\zs{0}\P(\xi>\alpha+z)\d z\\
&=\beta
\int^\infty_\zs{0}\P(\xi>\alpha+\beta x)\d x
\le
\beta
\int^\infty_\zs{0}\P(\Gamma^c_\zs{x})\d x
\le \beta M^*\,.
\end{align*}

\endproof

\subsection{ Proof of Theorem~\ref{Th.2.1}}\label{Se.A.1}

By making use of \eqref{1.1} and (\ref{2.3}) we obtain
\begin{equation}\label{A.1}
\|z-S\|^2=\gamma_n(z)\,+\,2\,F_n(z)\,+\,\|S\|^2\,,
\end{equation}
where $F_n(z)\,=\,n^{-1}\,\int^n_0\,z(t)\,\d \xi_t$.
Further, from the definition (\ref{2.7}), it follows that
for each $m\in\cM$
$$
\gamma_n(\wt{S}_\zs{\wt{m}})\,\le\,
\gamma_n(\wt{S}_\zs{m})\,+\,P_\zs{n}(m)-P_n(\wt{m})\,.
$$
Thus
\begin{equation}\label{A.2}
\|\wt{S}_\zs{\wt{m}}\,-\,S\|^2
\le
\|\wt{S}_m-S\|^2
+2F_n(\wt{z})+\varpi_\zs{n}(m,\wt{m})\,,
\end{equation}
where $\wt{z}=\wt{S}_\zs{\wt{m}} - \wt{S}_m$
and $\varpi_\zs{n}(m,\wt{m})= P_\zs{n}(m)-P_n(\wt{m})$.
Now for each  $x>0$, $0<\mu<1/2\lambda^*$ and set  $\iota\in\cM$
 we introduce the following gaussian function on $\cD_\zs{\iota}+\cD_\zs{m}$
\begin{equation}\label{A.3}
U_\zs{x,\iota}(z,\mu)=\frac{2F_n(z)}{\|z\|^2+\varrho^2_\zs{n,\iota}(x,\mu)}\,, \ \ \ \ \
z\in\cD_\zs{\iota}+\cD_\zs{m}\,,
\end{equation}
where
$$
\varrho_\zs{n,\iota}(x,\mu)\,=\,4\,
\sqrt{\frac{c(\mu)\,N+d_\zs{\iota}l_\zs{\iota}+x}{n \mu}}
\quad\mbox{with}\quad
c(\mu)=-\frac{1}{2}\,\ln(1-2\lambda^{*}\mu)
$$
and $N=\dim(\cD_\zs{\iota}+\cD_\zs{m})$.

Moreover, let functions $\phi_\zs{i_1},\ldots,\phi_\zs{i_N}$ be the subset of
$(\phi_j)_\zs{j\ge 1}$ which is a basis in
$\cD_\zs{\iota}+\cD_\zs{m}$. It should be noted that $N\le d_\zs{\iota}+d_\zs{m}$.
  Then
 one can write
 a normalized vector $\overline{z}=z/\|z\|$ for $z\ne 0$ as
$$
\overline{z}=\sum^N_{j=1}a_j\phi_\zs{i_j}
\quad\mbox{with}\quad
\sum^N_{j=1}\,a^2_j=1\,.
$$
Therefore
$$
U_\zs{x,\iota}(z,\mu)=\frac{2\, n^{-1/2}\,\|z\|}{\|z\|^2+
\varrho^2_\zs{n,\iota}(x,\mu)}\sum^N_{j=1}\,a_j\,\zeta_{i_j}
\quad\mbox{with}\quad
\zeta_j\,=\,\frac{1}{\sqrt{n}}\,\int^n_0\,\phi_\zs{j}\,\d \xi_t\,.
$$
By applying the Bunyakovskii-Cauchy-Schvartz
inequality one gets
\begin{equation}\label{A.4}
|U_\zs{x,\iota}(z,\mu)|\,\le\,
\frac{1}{\sqrt{n}\,\varrho_\zs{n,\iota}(x,\mu)}\,\eta_\zs{\iota}\,,
\end{equation}
where $\eta_\zs{\iota}=\sqrt{\sum^N_{j=1}\,\zeta^2_{i_j}}$.
Now note that by the condition $\C_1)$ the vector
$(\zeta_\zs{i_1},\ldots,\zeta_\zs{i_N})$ is  gaussian with
 zero mean and a non-generate covariance matrix $B_\zs{\iota}$.
Therefore
$$
\E\,e^{\mu\sum^N_{j=1}\zeta^2_{i_j}}
=\frac{1}{(2\pi)^{N/2}\sqrt{\det B_\zs{\iota}}}\,\int_{\R^N}\,
e^{-\frac{1}{2}x'T^{-1}_\zs{\iota}x}\,\d x\,,
$$
where $T_\zs{\iota}=(B^{-1}_\zs{\iota}-2\mu I_N)^{-1}$ and $I_\zs{N}$ is the identity matrix of order $N$.
One can easily verify that
\begin{align*}
\E\,e^{\mu\sum^N_{j=1}\zeta^2_{i_j}}=
\sqrt{\frac{\det T_\zs{\iota}}{\det B_\zs{\iota}}}=
\frac{1}{\sqrt{\det (I_N-2\mu\,B_\zs{\iota})}}\,.
\end{align*}
Thus, in view of the inequality
\begin{align*}
\det (I_N-2\mu B_\zs{\iota})\ge (1-2\mu\,\lambda_{max}(B_\zs{\iota}))^N
\end{align*}
and the condition $\C_2)$, we obtain 
\begin{align*}
\E\,e^{\mu\sum^N_{j=1}\zeta^2_{i_j}}\le (1-2\mu\lambda^*)^{-N/2}
=e^{c(\mu) N}\,,
\end{align*}
where $c(\mu)$ is defined in \eqref{A.3}.

Now, by the Chebyshev inequality,  for
 any $b>0$ and $0<\mu<1/2\lambda^*$, we obtain that
\begin{equation}\label{A.5}
\P(\,\eta_\zs{\iota}>b\,)\,
\le e^{c(\mu)N-\mu b^2}\,.
\end{equation}
Choosing in this inequality
$$
b=b_*(x,\iota)=\frac{1}{4}\sqrt{n}\,\varrho_\zs{n,\iota}(x,\mu)\,
=
\sqrt{\frac{c(\mu)\,N+d_\zs{\iota}l_\zs{\iota}+x}{\mu}}
$$
 yields
$$
\P(\eta_\zs{\iota}>b_*(x,\iota))\le\,e^{-x-d_\zs{\iota}l_\zs{\iota}}\,.
$$
Now let
$\Gamma_\zs{x}=\{\sup_{\iota\in\cM}\,
\eta_\zs{\iota}/b_*(x,\iota)\le 1\}$.
It is easy to see that
\begin{align*}
\P(\Gamma^c(x))\le
\sum_\zs{\iota\in\cM}\P(\eta_\zs{\iota}>b_*(x,\iota))
\le
 \sum_\zs{\iota\in\cM}e^{-x-d_\zs{\iota}l_\zs{\iota}}\,
= l^{*} e^{-x}\,.
\end{align*}
Thus, we obtain the following upper bound on the set $\Gamma_\zs{x}$
$$
\sup_{\iota\in\cM}\sup_{z\in\cD_\zs{\iota}+\cD_\zs{m}}
\,|U_\zs{x,\iota}(z,\mu)|\le 1/4\,,
$$
which implies
\begin{align*}
2\,F_n(\wt{z})&=
U_\zs{x,\wt{m}}(\wt{z},\mu)(\|\wt{z}\|^2\,+\,
\varrho^2_\zs{n,\wt{m}}(x))
\le \frac{1}{2}\|\wt{S}_\zs{\wt{m}}-S\|^2\,+\,
\frac{1}{2}\|\wt{S}_m-S\|^2\\
&+\frac{4}{n\mu}
(c(\mu) d_\zs{m}\,l_\zs{m}\,+\,
(c(\mu)\,+\,1)\,d_\zs{\wt{m}}l_\zs{\wt{m}})
+\frac{4}{n\mu}\,x\,.
\end{align*}
By making use of this inequality in (\ref{A.2}) we obtain,
on the set $\Gamma_\zs{x}$, that
\begin{align*}
\|\wt{S}_{\wt{m}}-S\|^2&\le \frac{1}{2}\,\|\wt{S}_{\wt{m}}-S\|^2
+\frac{3}{2}\|\wt{S}_m-S\|^2
+\varpi_\zs{n}(m,\wt{m})\\
&+\frac{4}{n\mu}
(c(\mu) d_\zs{m}\,l_\zs{m}
+(c(\mu) + 1)\,d_\zs{\wt{m}}l_\zs{\wt{m}})
+\frac{4}{n\mu}\,x\\
&= \frac{1}{2}\,\|\wt{S}_{\wt{m}}-S\|^2
\,+\,\frac{3}{2}\|\wt{S}_m-S\|^2+
\Omega(\rho,\mu)
+\frac{4}{n\mu}\,x\,,
\end{align*}
i.e.
$$
\|\wt{S}_{\wt{m}}-S\|^2
\le
3\|\wt{S}_m-S\|^2+
2\Omega(\rho,\mu)
+\frac{8}{n\mu}\,x\,,
$$
where
$$
\Omega(\rho,\mu)=\left(\rho+\frac{4c(\mu)}{\mu}\right)\frac{d_\zs{m}\,l_\zs{m}}{n}
+\left(\frac{4c(\mu)+4}{\mu}-\rho\right)\frac{d_\zs{\wt{m}}\,l_\zs{\wt{m}}}{n}\,.
$$
It is clear that to obtain a nonrandom minimal upper bound for this term
we have to resolve
the following optimazation problem
\begin{equation}\label{A.5-1}
\rho+\frac{4c(\mu)}{\mu}\to \min
\quad\mbox{subject to}
\quad
\frac{4c(\mu)+4}{\mu}-\rho\le 0\,.
\end{equation}
One can check directly that the solution of this problem
is $\mu=(z_\zs{*}-1)/2\lambda^{*}z_\zs{*}$
and the optimal value for $\rho$ is given in \eqref{2.6}.
 Thus, by choosing these parameters we have on the set $\Gamma_\zs{x}$ 
$$
\|\wt{S}_{\wt{m}}-S\|^2\le
3\|\wt{S}_m-S\|^2+16\lambda^{*}z_\zs{*}\frac{l_\zs{m}d_\zs{m}}{n}
+\frac{16z_\zs{*}\lambda^{*}}{n (z_\zs{*}-1)}\,x\,.
$$
Applying now Lemma~\ref{Le.A.1} 
with $\xi=\|\wt{S}_{\wt{m}}-S\|^2$,
$$
\alpha= 3\|\wt{S}_m-S\|^2+16\lambda^{*}z_\zs{*}l_\zs{m}d_\zs{m}/n\,,\quad
\beta =\frac{16z_\zs{*}\lambda^{*}}{n (z_\zs{*}-1)}
$$
and $M(x)=l^{*}e^{-x}$
  we obtain the inequality \eqref{2.9}.
Hence Theorem~\ref{Th.2.1}.
\endproof

\subsection{Some properties of the Fourier coefficients}\label{Se.A.4}

\begin{lemma}\label{Le.A.2}
Let $S$ be a function in $\C^{k}[0,1]$
such that $S^{(j)}(0)=S^{(j)}(1)$ for all $0\le j\le k$ and, such that,
for some contants $L_0>0$, $L>0$ and $0\le\,\alpha\,<1$
\begin{equation}\label{A.9}
\max_{0\le j\le k-1}\max_{0\le x\le 1}|S^{(j)}(x)|\le L_0 \ \ \mbox{and} \ \
|S^{(k)}(x)-S^{(k)}(y)|\le L\,|x-y|^\alpha
\end{equation}
for all $x\,,y\in [0,1]$.
Then the Fourier coefficients $(a_k)_\zs{k\ge 0}$ and $(b_k)_\zs{k\ge 1}$
 of the function $S$, defined as
$$
S(x)=\frac{a_0}{2}+\sum^\infty_{k=1}
(a_k\cos(2\pi\,kx)+b_k\sin(2\pi\,kx))\,,
$$
 satisfy the following inequality
\begin{equation}\label{A.10}
\sup_{n\ge 0}
(n+1)^\beta\left(\sum^\infty_{j=n}
\,(a^{2}_j+b^{2}_j)\right)^{1/2}\le c^*\, (L+L_0)\,,
\end{equation}
where $\beta=k+\alpha$ ($k$ being an integer and $0\le\alpha<1$) and
$$
c^*=1+2^\beta+\pi^4\,9^\beta\,
\frac{\int^\infty_0\,u^{\alpha-3}
\,\sin^4(\pi\,u)\,\d u}
{8\,\int^{1/2}_0\,
u^{-4}\,\sin^4(\pi\,u)\,\d u}\,.
$$
\end{lemma}
\noindent Proof of this Lemma is given in Fourdrinier and pergamenshchikov (2007) Appendix A4.

\medskip

\subsection{Lower bound for the parametric  model }\label{Se.A.4-1}
We consider in this section the following model
\begin{equation}\label{A.10-1}
\d y_\zs{t}=S(t,z)\d t+\d w_\zs{t}\,,
\end{equation}
where $(w_\zs{t})_\zs{t\ge 0}$ is a standard brownian motion;
 $z\in\bbr^m$ is unknown vector parameter. Let now
$\pi$ be a prior distribution density on $\bbr^m$ of the form
$$
\pi(x)=\prod^m_\zs{l=1}\,\pi_\zs{l}(x_\zs{l})\,,
$$
where $\pi_\zs{l}$ is a positive density on the interval
$[-\delta_\zs{l},\delta_\zs{l}]$ for some $\delta_\zs{l}>0$.
This means that the density $\pi$ has the following support
$$
\Gamma=[-\delta_\zs{1},\delta_\zs{1}]\times\ldots\times
[-\delta_\zs{m},\delta_\zs{m}]\,.
$$
We set
\begin{equation}\label{A.10-2}
\zeta_\zs{l}(z)=\int^n_\zs{0}\,\frac{\partial}{\partial
z_\zs{l}}\,S(t,z)(\d y_\zs{t}-S(t,z)\d t)\,.
\end{equation}
Now we give some version of the van Trees inequality Gill and Levit (1995)
for process \eqref{A.10-1}.
\begin{lemma}\label{Le.A.2-1}
For each $l\ge 1$, any estimator $\wt{z}_\zs{l}$ based on
observations $(y_\zs{t})_\zs{0\le t\le n}$
 satisfies the
following inequality
\begin{equation}\label{A.10-3}
\int_\zs{\Gamma}\,\E_\zs{S_\zs{z}}(\wt{z}_\zs{l}-z_\zs{l})^2\,
\pi(\d z)\,\ge\,\frac{1}{A_\zs{l}+B_\zs{l}}\,,
\end{equation}
where $\E_\zs{S_\zs{z}}$ denotes the expectation with respect to the
distribution of process \eqref{A.10-1},
$$
A_\zs{l}
=\int_\zs{\Gamma}\E_\zs{S_\zs{z}}\,\zeta^2_\zs{l}(z)\,\pi(z)\d z
\quad\mbox{and}\quad
B_\zs{l}
=\int^{\delta_\zs{l}}_\zs{-\delta_\zs{l}}\,
\frac{\left(\dot{\pi}_\zs{l}(u)\right)^2}{\pi_\zs{l}(u)}\,\d u\,.
$$
\end{lemma}
\proof
It will be noted that the density of
the distribution of process \eqref{A.10-1}
with respect to the Wiener measure $\mu_\zs{w}$ 
on $\cY=\C[0,n]$ is defined as
$$
f(y,z)=e^{\int^n_0\,S_z(t)\,\d y_t\,
-\,\frac{1}{2}\,\int^n_0\,S^2_z(t)\,\d t}\,.
$$
Therefore, by applying the method from Gill and Levit (1995) we obtain
the lower bound \eqref{A.10-3}.
 Hence Lemma~\ref{Le.A.2-1}.
\endproof

\subsection{Proof of Theorem~\ref{Th.6.1}}\label{Se.A.5}

To prove this theorem we adapt the proof of Theorem~\ref{Th.2.1} for this case.
In this case equality \eqref{A.1} becomes
$$
\|z-S\|^2_p=\gamma_\zs{n,p}(z)\,+\,2\,F_\zs{n,p}(z)\,+\,2\,G_\zs{p}(z,S)\,
+\,\|S\|^2_p\,,
$$
with
\begin{align*}
F_\zs{n,p}(z)\,=\,\frac{1}{n}\,\sum^{np}_\zs{k=1}\,z(t_k)\,\Delta \xi_\zs{t_k}\,,\ \
G_\zs{p}(z,S)\,=\,\frac{1}{p}\sum^{p}_\zs{k=1}\,z(t_k)\,h_\zs{k}(S)\,,
\end{align*}
where the sequence $h_\zs{k}(S)$ is defined in \eqref{6.10}.
Similarly to the proof of Theorem~\ref{Th.2.1}, one can show that
\begin{align}\nonumber
\E_\zs{S}\|\wt{S}_\zs{\wt{m}_\zs{p},p}-S\|^2_p\,&\le\,
3\E_\zs{S}\|\wt{S}_{m,p}-S\|^2_p
+16\lambda^{*}z_\zs{*}\frac{d_ml_\zs{m}}{n} \\ \label{A.11}
&+4\E_\zs{S}|G_\zs{p}(\wt{z}_p,S)|+
\frac{\lambda^{*}\tau_\zs{0}}{n}\,,
\end{align}
where $\wt{z}_p=\wt{S}_\zs{\wt{m}_\zs{p},p}-\wt{S}_{m,p}$.
Now we note that for any $\nu>0$
\begin{align*}
2G_p(\wt{z},S)\,&\le \nu \|\wt{z}\|^2_p+\nu^{-1}\,H_p(S)\\
&\le 2\nu\|\wt{S}_\zs{\wt{m}_\zs{p},p}-S\|^2_p+2\nu\|\wt{S}_{m,p}-S\|^2_p
+\nu^{-1}\,H_p(S)\,.
\end{align*}
Therefore, taking into account the last inequality in
\eqref{A.11},
 we obtain the following upper bound
\begin{align*}
\E_\zs{S}\|\wt{S}_\zs{\wt{m}_\zs{p},p}-S\|^2_p\,&\le\,
\frac{3+4\nu}{1-4\nu}
\,\|\wt{S}_\zs{m,p}-S\|^2_p+\frac{16\lambda^{*}z_\zs{*}}{1-4\nu}\,\frac{d_m l_\zs{m}}{n}
\\
&+\frac{\lambda^{*}\tau_\zs{0}}{(1-4\nu)n}+\frac{2}{\nu(1-4\nu)}\,H_p(S)\,.
\end{align*}
By minimzing the last term with respect to $\nu$ (i.e. maximazing $\nu(1-4\nu)$)
 we find that $\nu=1/8$. Thus the last inequality implies
the upper bound \eqref{6.10}.
Hence, Theorem~\ref{Th.6.1}.
\endproof

\subsection{Proof of \eqref{6.15}}\label{Se.A.6}
Consider first the principal term in \eqref{6.13}.
Let $(s_j)_\zs{j\ge 1}$ be the Fourier coefficients for $S$ in
$\cL_2[0,1]$ used in $\Theta_\zs{\beta,r}$.
By setting 
$\Delta_{j}(t)=S-\sum^j_\zs{i=1} s_i \phi_i(t)$
, one can estimate
$\|S_\zs{m_j,p}-S\|^2_p$ as
\begin{align*}
\|S_\zs{m_j,p}-S\|^2_p\,=\,\inf_\zs{a_1,\ldots,a_j}\,
\|S-\sum^j_\zs{i=1}\,a_j\,\phi_i\|^2_p\,\le\,\|\Delta_j\|^2_p\,.
\end{align*}
By the definition of
$\varsigma_j(S)$ in (\ref{4.1-1}) we obtain that
\begin{align}\nonumber
\|\Delta_j\|^2_p\,
&\le\,2\,\int^1_\zs{0}\,\Delta^2_j(t)\,\d t\,+\,2\,
\sum^p_\zs{k=1}\,\int^{t_k}_\zs{t_\zs{k-1}}\,(\Delta_j(t_k)- \Delta_j(t))^2\,\d t\\ \label{A.11-1}
&=\,2\,\varsigma_\zs{j+1}(S)\,+\,2\sum^p_\zs{k=1}\,\int^{t_k}_\zs{t_\zs{k-1}}\,
\left(\int^{t_k}_\zs{t}\dot{\Delta}_j(u)\,\d u\right)^2\,\d t\,.
\end{align}
Moreover, the Bunyakovskii-Cauchy-Schwartz inequality implies that
$$
\|\Delta_j\|^2_p\le 2\,\varsigma_\zs{j+1}(S)\,+\,\frac{2}{p^2}\,\|\dot{\Delta}_j\|^2\,.
$$
Notice now that for the trigonometric basis \eqref{5.1}
and for the functions $S$ from $\Theta_\zs{\beta,r}$ with $\beta>1$
we obtain that for any $j\ge 0$
\begin{align}\nonumber
\|\dot{\Delta}_j\|^2&=
\sum^{\infty}_\zs{i=j+1}\,s^2_i\,\|\dot{\phi}_i\|^2
\le
\pi^2\,\sum^{\infty}_\zs{i=j+1} s^2_i i^2\\ \nonumber
&\le
\pi^2 (j+1)^2\,\varsigma_\zs{j+1}(S) + 2 \pi^2\,\sum_\zs{i>j}\,(i+1)\,
\varsigma_\zs{i+1}(S)\\ \label{A.11-2}
&\le\,r^2\,
\pi^2\,(j+1)^{-2(\beta-1)}\,\frac{\beta}{\beta-1}\,.
\end{align}
Therefore for $1\le j\le p$
\begin{align}\label{A.12}
\sup_\zs{S\in \Theta_\zs{\beta,r}}\,
\|S_\zs{m_j,p}-S\|^2_p\,&\le\,2\,\frac{r^2}{j^{2\beta}}\,
\left(1+\,\frac{\pi^2 j^2}{p^2}\,\frac{\beta}{\beta-1}\right)\,.
\end{align}
Moreover, taking into account that
$H_\zs{p}(S)\le\,p^{-2} \|\dot{S}\|^2$, through \eqref{A.11-2}
with $j=0$ we get that for $p\ge n^{1/2}$
$$
\sup_\zs{S\in \Theta_\zs{\beta,r}}\,H_p(S)
\le \frac{\beta \pi^2\,r^2}{p^2(\beta-1)}
\le \frac{\beta \pi^2\,r^2}{n(\beta-1)}\,.
$$
Thus \eqref{6.13} implies  that
for any $\varepsilon>0$
there exists some constant $C^*=C^*(r,\epsilon)>0$ such that for any
$\beta\ge 1+\epsilon$, $p\ge n^{1/2}$ and for  $1\le j\le p$
\begin{align*}
\cR_\zs{n,p}(\wh{S}_\zs{\wh{m}_\zs{p},p},\beta)\,
\le\,C^*\,\left(\,\left(j^{-2\beta}\,+\,\frac{j}{n}\right)\,n^{\frac{2\beta}{2\beta+1}}
\,+\,n^{-\frac{1}{2\beta+1}}
\right)\,.
\end{align*}
This bound with
$j=j_*=[n^{\frac{1}{2\beta+1}}]+1$ implies immediately
inequality \eqref{6.15}.
\endproof

\subsection{Proof of \eqref{6.16}}\label{Se.A.7}

Notice now that
for any estimator $\wt{S}_n$ by putting
$T_p(\wt{S}_n)(t)=\sum^p_\zs{j=1}\wt{S}_n(t_j)
\Chi_\zs{(t_\zs{j-1},t_j]}$ we can represent
the accuracy of this estimator as
$$
\|\wt{S}_n-S\|^2_\zs{p}=\sum^p_\zs{j=1}\,
\int^{t_j}_\zs{t_\zs{j-1}}(T_p(\wt{S}_n)(t)\,-S(t_k))^2\,\d t\,.
$$
Therefore, for any $0<\epsilon<1$ we can estimate with below this accuracy as

\begin{align*}
\|\wt{S}_n-S\|^2_\zs{p}
\ge
(1-\epsilon)\|T_p(\wt{S}_n)-S\|^2
-(\epsilon^{-1}-1) \sum^p_\zs{j=1}
\int^{t_j}_\zs{t_\zs{j-1}} (S(t)-S(t_j))^2\d t\,.
\end{align*}
Moreover, similarly to \eqref{A.11-1}--\eqref{A.12}
 we obtain that
$$
\sum^p_\zs{j=1}\,
\int^{t_j}_\zs{t_\zs{j-1}}\,(S(t)-S(t_j))^2\d t\,
\le\,p^{-2}\,\|\dot{S}\|^2\,\le\,r\pi^2 n^{-1}\,.
$$
Therefore
$$
\cR_\zs{n,p}(\wh{S}_n,\beta)\,\ge\,(1-\epsilon)\,
\inf_\zs{T_\zs{n}}\,\cR_\zs{n}(T_\zs{n},\beta)\,-\,(\epsilon^{-1}-1)\,
r\pi^2\,n^{-\frac{1}{2\beta+1}}\,,
$$
where the risk $\cR_\zs{n}(T_\zs{n},\beta)$ is defined by \eqref{4.2} for some estimator
$T_\zs{n}$.
Now Theorem~\ref{Th.5.1} directly implies \eqref{6.16}.
\endproof


\end{document}